\documentclass[aip,graphicx]{revtex4-1}

\draft 

\usepackage{graphicx}
\usepackage{amsmath,amssymb}
\usepackage[version=4]{mhchem}
\usepackage{siunitx}
\usepackage{multirow}
\usepackage[utf8]{inputenc}
\usepackage[english]{babel}

\def\krmi{k_a}
\def\krme{\krmi^{\rm meso}}

\def\rrad{\sigma}

\def\microreact{\tau_{\rm{react}}^{\rm{micro}}}
\def\mesoreact{\tau_{\rm{react}}^{\rm{meso}}}

\begin{document}


\title{Hierarchical Reaction-Diffusion Master Equation} 



\author{Stefan Hellander}
\email[]{stefan.hellander@it.uu.se}
\affiliation{Department of Information Technology, Uppsala University, Box 337, SE-755 01, Uppsala, Sweden}
\author{Andreas Hellander}
\email[]{andreas.hellander@it.uu.se}
\affiliation{Department of Information Technology, Uppsala University, Box 337, SE-755 01, Uppsala, Sweden}


\date{\today}
\begin{abstract}

We have developed an algorithm coupling mesoscopic simulations on different levels in a hierarchy of Cartesian meshes. Based on the multiscale nature of the chemical reactions, some molecules in the system will live on a fine-grained mesh, while others live on a coarse-grained mesh. By allowing molecules to transfer from the fine levels to the coarse levels when appropriate, we show that we can save up to three orders of magnitude of computational time compared to microscopic simulations or highly resolved mesoscopic simulations, without losing significant accuracy. We demonstrate this in several numerical examples with systems that cannot be accurately simulated with a coarse-grained mesoscopic model. 

\end{abstract}

\pacs{}
\maketitle 

\section{Introduction}

Spatial stochastic modeling of reaction-diffusion kinetics is a popular tool to study the fine-grained molecular details of intracellular regulation. By being able to capture both the inherent spatial aspects of signal transduction as well as the discrete and stochastic nature of molecular interaction in the low-copy number regime, these types of models offer the possibility of insights not attainable by either more detailed molecular dynamics models (due to their computational cost) or phenomenological macroscopic models  (due to the deterministic description)\cite{sturrockSpatialStochasticModelling2013,sturrockRoleDimerisationNuclear2013,vanzonDiffusionTranscriptionFactors2006,elfSpontaneousSeparationBistable2004,fangeNoiseInducedMin2005,lawsonSpatialStochasticDynamics2013}.

Spatial stochastic simulation algoritms used in systems biology have to carefully balance the need for high spatial resolution with the need for a low computational cost, in order to study regulatory processes over physiologically relevant time scales (entire cell cycles). Two model formalisms have attracted particular attention in the field: the mesoscopic on-lattice Reaction-Diffusion Master Equation (RDME) and the microscopic off-lattice Collins-Kimball-Smoluchowski (CKS) model.  In the former, proteins are modeled as point particles and are diffusing on the grid according to a discrete jump process, and they are able to react when finding themselves in the same voxel, whereas in the latter, proteins are modeled as individual hard spheres and diffuse continuously in space according to Brownian motion. The CKS model is generally considered being a more accurate model than the RDME, although there is no formal relationship between these two models in the sense that one arises as an approximation of the other. 

By choosing reaction rates in the RDME so that the properties of the microscopic model are captured, it is possible to relate the two models to each other formally \cite{hellanderReactionRatesMesoscopic2015,hellanderReactiondiffusionMasterEquation2012,fangeStochasticReactiondiffusionKinetics2010}. With the choice of mesoscopic rate constants from \cite{hellanderReactionRatesMesoscopic2015} it is possible to match the mean binding time between two molecules in the two models, down to a critical size of the mesh \cite{hellanderReactiondiffusionMasterEquation2012}. For meshes finer than this critical size, the on-lattice RDME cannot capture the microscopic dynamics accurately. 

Due to the popularity of these modeling frameworks, several capable open-source software frameworks have been developed to support spatial stochastic modeling both for the RDME \cite{drawertMOLNsCloudPlatform2015,drawertURDMEModularFramework2012,hattneStochasticReactiondiffusionSimulation2005} and the particle-based model \cite{andrewsDetailedSimulationsCell2010,kerrFASTMONTECARLO2008,PhysRevLett.97.230602}. Some software support simulation on multiple levels and integrate one or more spatial stochastic simulators, such as VCell \cite{schaffGeneralComputationalFramework1997} and StochSS \cite{drawertStochasticSimulationService2016}. 

In summary, several capable tools exist to simulate a reaction-diffusion system either on the microscopic or mesocopic scale. The mesoscopic on-lattice model offers superior simulation speed, assuming sufficient accuracy can be obtained with a relatively coarse mesh. However, studies have highlighted scenarios where a very high spatial resolution is necessary to capture microscopic properties such as the rebinding time distributions accurately, in order to capture the correct macroscopic behavior  \cite{vanzonDiffusionTranscriptionFactors2006, Mahmutovic2012,fangeStochasticReactiondiffusionKinetics2010}. For these systems, the computational cost becomes substantial both for microscopic particle-based methods and for on-lattice simulations with high spatial resolution \cite{hellanderReactionRatesMesoscopic2015}. 

The reason for the rapid growth in computational cost differs in simulation based on the CKS and RDME models. Implementations of the microscopic model, being a many-body problem,  scales poorly with the number of particles in the simulation.  The Greens Function Reaction Dynamics (GFRD) algorithm improves performance over naive Brownian dynamics for sparse systems with relatively few particles \cite{vanzonGreenSfunctionReaction2005}, but the computational cost still becomes overwhelming for systems with many interacting particles. If the mesh used in the RDME can be chosen relatively coarse, simulations on the mesoscopic scale is typically orders of magnitudes faster than simulations based on the CKS model and scales linearly with the number of particles. However, the on-lattice RDME suffers from stiffness, leading to a quickly growing computational cost as the mesh is refined. The number of diffusive jumps per time unit of simulation time is proportional to $D/h^2$, where $D$ is the diffusion constant and $h$ is the length of a voxel \cite{hellanderLocalErrorEstimates2014}. 

The problem for practical modeling is that spatial models often have at least a few reactions that are diffusion limited and hence require a high spatial resolution, but on the same time species that are present in relatively large copy numbers (in the hundreds or thousands). This results in a situation where neither method performs well. A natural way to approach these multiscale systems are to blend mesoscopic and microcopic methods in one single simulation. Previous work on such hybrid methods have highlighted the large computational savings made possible by a multiresolution approach \cite{hellanderMesoscopicmicroscopicSpatialStochastic2017}. A challenge is to partition the system into its microscopic and mesoscopic part without prior knowledge about the system dynamics. In previous work Hellander et al. demonstrated how analysis from \cite{hellanderLocalErrorEstimates2014} can be used for automatic system partitioning \cite{hellanderMesoscopicmicroscopicSpatialStochastic2017}.

In addition to speeding up simulations with multiscale reaction properties, mesoscopic-microscopic hybrid methods can be applied to split simulation accuracy in different parts of the domain \cite{hellanderCoupledMesoscopicMicroscopic2012,hybrid2,tworegime1}, and to augment mesoscopic models for situations where the mesocopic framework is not well defined, such as for interactions between molecules and surfaces, and for 2D-3D interactions \cite{hellanderCoupledMesoscopicMicroscopic2012}.   

Hybrid methods can achieve good speedups,  but a distinct disadvantage is the relatively large complexity in their implementation, and overhead caused by switching between data structures optimal for the respective algorithms. In this paper we present a new pure on-lattice multiscale and multilevel method for spatial stochastic simulations. Based on our previous analysis of the accuracy of the the RDME on different spatial resolutions \cite{hellanderReactionRatesMesoscopic2015,hellanderReactiondiffusionMasterEquation2012}
we design a hierarchical simulation algorithm that employs several meshes of different resolution in order to capture the fine scale dynamics of highly diffusion-limited reactions while avoiding the need to resolve the entire systems on that same high level. In a series of numerical examples of increasing complexity, we demonstrate an accuracy comparable with pure GFRD simulations at a simulation cost up to three orders of magnitude below state-of-the-art GFRD implementations.

\section{Background}

In the next section, we describe a method that allows reactions to take place on different mesh resolutions depending on the degree of diffusion control. This hierarchical RDME model allows for high accuracy at a much reduced cost compared to a fully microscopically resolved system, for models with multiscale properties. In this section, we first describe the underlying mesoscopic model, and then briefly review the microscale model. We consider the more fine-grained microscale as the correct model later when computing the error of the mesoscopic simulations.

\subsection{Reaction-diffusion master equation}

The reaction-diffusion master equation is the natural spatial extension of the popular well-mixed Markov process description of chemical kinetics \cite{vankampenStochasticProcessesPhysics1992, gillespiePerspectiveStochasticAlgorithms2013}. This model formalism is widely used in systems biology, and models the state $\mathbf{x}$ of the system as a vector consisting of the discrete number of molecules of each chemical species. Formally, chemical species $X_i, i=1\ldots N$ participate in M chemical reactions $R_{j}, j=1\ldots M$. For example, a bimolecular reaction where species $X_1$ react with $X_2$ to form $X_3$ can be written as 
\begin{align}
X_{1} + X_{2} \xrightarrow {k} X_3
\label{eq:bimolecular}
\end{align}
\noindent 
where $k$ is the reaction rate parameter for the reaction. Using mass action kinetics, the propensity function for the reaction \eqref{eq:bimolecular} is a function of the rate constant and the copy number of the reactants $X_1, X_2$, $a(\mathbf{x}) = kx_1x_2$. In the Markov process formalism, the inverse of the propensity $1/a_r(\mathbf{x})$ gives the transition rate for changing states from $[x_1,x_2,x_3]$ to $[x_1-1,x_2-1, x_3+1]$.     

The time evolution of the probability density of the system is governed by the forward Kolmogorov equation, or the chemical master equation (CME), but since this equation is infeasible to solve for systems with a large number of chemical species, kinetic Monte Carlo simulation using the direct stochastic simulation algorithm (SSA) \cite{gillespieGeneralMethodNumerically1976}, or one of its many optimized or approximate variants \cite{gillespiePerspectiveStochasticAlgorithms2013}, is normally used to analyze the system.  

In the spatial stochastic case, the computational domain is partitioned into K voxels $V_k$ using a mesh. Molecules move by diffusion, modeled as discrete jump events between adjacent voxels, according to a linear event
\begin{align}
X_{ij} \xrightarrow{d_{ijk}} X_{ik}
\end{align}
The rate $d_{ijk}$ depends on the diffusion constant of $X_i$ and on the shape and size of the voxels \cite{engblomSimulationStochasticReactionDiffusion2009}. Chemical reactions are modeled as in the well-mixed model, but now locally confined to individual voxels. Compared to simulation of well-mixed systems, the computational cost grows quickly with the size of the mesh. If $h$ is a measure of the length scale of the voxel, the total number of diffusion events in a simulation scale like $1/h^2$. This stiffness problem causes RDME simulations to become highly computationally expensive if a high spatial resolution is needed.

\subsection{Next-Particle Method}

There are different methods for generating trajectories of the RDME. Widely used is the Next-Subvolume Method (NSM), in which the population count of each species is tracked inside the voxels \cite{elfSpontaneousSeparationBistable2004}. While the specific choice of solver is not critical for the hierarchical RDME method (hRDME), we here choose a different approach for practical reasons. It will be useful to know for how long each molecule has existed within the system, and therefore we implement the hRDME with the Next-Particle Method (NPM) \cite{hellanderMesoscopicmicroscopicSpatialStochastic2017} as the underlying mesoscopic solver. In this section we describe the method as implemented on a single mesh, and then, in Section \ref{sec:method}, we describe how it can be adapted to the case of multiple meshes.

The method is initialized by sampling a voxel for each molecule in the system. Usually the initial distribution of molecules is uniform, in which case we sample a voxel from a uniform distribution, but this is not a necessary requirement. Just as in the case of the NSM we will maintain an event queue sorted in descending order based on the time for each event (so that the next event in the queue is also the event that fires next). For each molecule in the system we therefore sample tentative events and add them to the queue.

\subsubsection{Initialization}

\begin{enumerate}
	\item For each molecule, add a tentative next diffusion event.
	\item For all molecules participating in a unimolecular event, add the tentative next unimolecular event. Note that we only add one tentative event for each molecule, as later unimolecular events could never fire.
	\item For all molecules participating in a bimolecular event, add the tentative next bimolecular event. 
	\item In addition, we add tentative events of the type $\emptyset \to S$. For each reaction of that type we add one tentative next reaction.
\end{enumerate}

\subsubsection{Propagation}

After the system has been initialized, the algorithm proceeds by executing the events in order.

\begin{enumerate}
	\item If the next event is a diffusion event, move the molecule accordingly. If the molecule participates in any bimolecular reactions, remove those from the queue. Unimolecular events are not affected by diffusion events, so they are left on the queue. Finally, sample new tentative bimolecular events and a new diffusion event.
	\item If the next event is a unimolecular event, remove all tentative events involving the reacting molecule. Initialize new molecules inside the same voxel as the reacting molecule, and finally sample new tentative diffusion and reaction events for each new molecule.
	\item If the next event is a bimolecular event, remove all tentative events involving either of the reacting molecules. Initialize products inside the same voxel, and sample new tentative events for all the products.
	\item If the next event is of the type $\emptyset\to S$, initialize a new molecule of species $S$ into the system. Unless otherwise specified, its initial position will be uniform.
	\item Repeat until the final time $T$.
\end{enumerate}

\subsection{Microscopic scale}

On the mesoscopic scale, particles are restricted to nodes on a computational grid. In contrast, on the microscopic scale, particles diffuse freely in continuous space according to normal diffusion. On the mesoscopic scale, particles are point particles, while on the microscopic scale they are modeled by hard spheres. Here reactions occur with some probability when the molecules collide. The reaction dynamics is governed by the probability density function solving the Collins-Kimball-Smoluchoski PDE \cite{Smol,CollinsKimball,vanzonGreenSfunctionReaction2005}.

Let $r$ be the distance between two reactive molecules $A$ and $B$, $D$ the sum of the diffusion constants, $\sigma$ the sum of the reaction radii, and $k_r$ the reaction rate. The probability for the distance $r$ at time $t$, given that the distance was $r_0$ at $t_0$, is given by $p(r,t|r_0,t_0)$, solving the equation
\begin{equation}
\frac{\partial p}{\partial t} = D\Delta p(r,t|r_0,t_0)
\end{equation}
with boundary condition
\begin{equation}
K\frac{\partial p}{\partial n}\bigg|_{r=\sigma} = k_r p(r,t|r_0,t_0),
\end{equation}
where
\begin{align}
K = \begin{cases}
4\pi\sigma^2 D\,\, \mathrm{(3D)}\\
2\pi\sigma D\,\, \mathrm{(2D)}.
\end{cases}
\end{align}

There exist several popular implementations of solvers of this model. Prominent examples are Smoldyn \cite{smoldyn}, MCell \cite{kerrFASTMONTECARLO2008}, and eGFRD \cite{vanzonGreenSfunctionReaction2005}. The former two, Smoldyn and MCell, take a similar approach in that they select a fixed time step, and proceed by propagating the system one time step at a time. An alternative approach is implemented in eGFRD where the system is propagated in continuous time. This approach tends to be more efficient if very high accuracy is required and if the system is reasonably sparse, while Smoldyn and MCell can be significantly more efficient in other cases.

In this paper we determine the accuracy of our simulations by comparing to corresponding simulations on the microscopic scale. To ensure the highest possible accuracy on the microscopic scale, we have compared to results obtained with the eGFRD algorithm, as well as another efficient implementation of a similar algorithm \cite{HELLANDER20113948,hellanderSingleMoleculeSimulations2013}.

\section{Method}
\label{sec:method}

We are often interested in simulating systems displaying dynamics on widely different scales. Parts of the system require a high spatial resolution, while other parts can be simulated on a coarse-grained mesh to satisfactory accuracy.

Instead of simulating the whole system on the fine-grained level, we will here describe an approach to coupling several mesh resolutions. Some molecules will be simulated on a fine-grained mesh, while others can be simulated to high accuracy on a much more coarse-grained mesh. Molecules can also be initialized on a fine-grained mesh, and after diffusing for a sufficiently long time, be transfered to a more coarse-grained mesh.

\subsection{A hierarchy of meshes}

Instead of simulating the entire system on a single mesh, we introduce a hierarchy of meshes. Depending on the dynamics of the system, some molecules may require a very high spatial resolution, while we can get away with simulating others on a much coarser mesh.

Here, for simplicity, we will consider Cartesian meshes only. The coarsest possible mesh is a single voxel. This corresponds to a fully well-mixed system. This mesh can then be successively refined by halving the voxel width, thus obtaining a sequence of meshes with 1, $2^3$, $4^3$, $8^3$, \ldots, number of voxels in 3D. By halving the width of the voxels in each step, each voxel will be fully contained within a voxel on a coarser mesh. This is not a neccessary requirement for the method to work, but it does simplify the implementation and keeps the overhead of the method at a minimum.

With this particular structure of the hierarchy of meshes, in which each voxel on a finer mesh is fully contained within a voxel on a coarser mesh, it is fairly straightforward to map molecules between the different scales. The mapping is a pure preprocessing step, in which each voxel on each mesh is assigned a parent voxel in the mesh one level coarser, and children voxels in the mesh one level finer. This assignment is particularly simple for Cartesian meshes, but would be possible to perform also in the case of an unstructured mesh (albeit much more computationally expensive).

\subsection{Move molecules between meshes}

The core idea of the algorithm is to transfer molecules between the different levels of the hierarchy depending on the dynamics of the system. For each species we can determine the finest mesh resolution necessary to resolve all dynamics involving that species (see Section \ref{choosehier}), and each molecule of that species will be initalized on the that mesh resolution. Depending on how the simulation proceeds, the molecule can be moved to coarser levels in the hierarchy, and products resulting from reactions involving that molecule can be moved to finer levels in the hierarchy.

\subsubsection{Move a molecule from a fine mesh to a coarse mesh}
However, if the molecule survives for long enough, and thus diffuse enough, it can successively be moved to a more coarse-grained mesh resolution without losing too much accuracy. The time until we can move a molecule from a fine mesh to a coarser mesh is related to the diffusion constant $D$ of the molecule and the width $h$ of the voxels in the current mesh. In particular, the time $t_{\rm{transfer}}$ until we can transfer a molecule to a coarser level is given by the relation
\begin{equation}
\label{Ceq}
h = \sqrt{\frac{6Dt_{\rm{transfer}}}{C}}
\end{equation}
for some constant $C$. In words, the molecule should, on average, diffuse a distance that is a multiple $\sqrt{C}$ of the voxel width on the current mesh, before we move it to a coarser mesh. This corresponds to the molecule getting ``well-mixed'' on the length scale of the voxels. 

A voxel on the fine mesh is fully contained within a voxel on the coarse mesh. When a molecule is transfered to a coarser mesh, it is simply placed in the voxel containing its current voxel. In Section \ref{sec:rebindingdyn} we determine a reasonable value for $C$ numerically.

\subsubsection{Move a molecule from a coarse mesh to a fine mesh}

Sometimes a molecule needs to be moved from a coarse mesh to a finer mesh in the hierarchy. For instance, when a molecule dissociates on a coarse mesh, the products might have to be initialized on a much finer mesh than that occupied by the reacting molecule. This is done by placing the molecule randomly inside one of the voxels contained within the voxel on the coarse mesh. 

\subsubsection{When can a molecule be moved?}

A molecule is only transfered between meshes immediately after a diffusion event has fired, and before we sample new tentative bimolecular reactions and a new tentative diffusion event.

The reason is that transfering a molecule to a new mesh is similar to executing a diffusion event; the molecule is placed inside a new voxel. This means that after transfering a molecule, we have to sample a new tentative bimolecular event and a new diffusion event. Transfering the molecule immediately following a diffusion event thus minimizes the overhead, because we need to perform these operations either way. We also avoid introducing a bias by artifically discarding tentative reaction events following a molecule transfer.

\subsection{Reactions}

When simulating the RDME on a single mesh, bimolecular reactions may occur when molecules occupy the same voxel, and products of zeroth- and first-order reactions are simply placed in the voxel of the reacting molecule. In the case of the hRDME, it will not be as straightforward. In particular, we need to determine the reaction rate for two reactive molecules occupying overlapping voxels on different levels in the hierarchy of meshes. Also, molecules can be initialized on a different level than that occupied by the reacting molecules, and in those cases we need to determine which voxel the products should be placed inside.

\subsubsection{Zeroth order reactions}

Reactions of the type $\emptyset\to S$ are executed just as in the standard NSM algorithm. We sample a tentative reaction time, and add the tentative event to the reaction queue. If the reaction fires, the new molecule is initialized on the mesh size required for the species $S$.

\subsubsection{Diffusion events}

A diffusion event is executed by first moving the molecule to one of the neighboring voxels with a uniform probability. Following a diffusion event, all tentative bimolecular events involving the molecule is removed from the queue. If we find a new tentative bimolecular event in the updated voxel, it is added to the queue. Unimolecular reactions are not affected by the diffusion event. A new tentative diffusion is also added to the queue.

\subsubsection{Unimolecular reactions}

Whenever a molecule is introduced into the system, we sample the next tentative unimolecular reaction involving this molecule. We will not have to update this event during the lifespan of the molecule, as unimolecular reactions are not affected by diffusion events.

When a unimolecular event fires, the reacting molecule is replaced by the product molecules. With each product molecule is associated a required mesh size. The products are initialized on their respective required mesh size. If this mesh is finer than the mesh occupied by the reacting molecule, we sample a voxel uniformly from the children of this voxel. If the mesh is coarser than that of the reacting molecule, we sample a voxel from the parents of the current voxel occupied by the reacting molecule.

For each new molecule introduced, we sample the corresponding tentative next events (diffusion, unimolecular, and bimolecular).

\subsubsection{Bimolecular reactions}

Let $S_1$ and $S_2$ be two species that react according to $\ce{S_1 + S_2 ->[k_a] S_3}$. If two molecules of species $S_1$ and $S_2$ occupy the same voxel on the same level in the hierarchy, they react just as in the NPM. However, here we frequently encounter the case where one of the molecules, $M_1$ of species $S_1$, occupies a voxel on a finer mesh than that occupied by molecule $M_2$ of species $S_2$. There is still a possibility that the molecules react if the voxel occupied by $M_1$ is contained within the voxel occupied by $M_2$. In this case they will react as if both molecules occupied the voxel occupied by $M_2$ on the coarser mesh.

Assume that the molecules react with the rate $k_r^{\rm{meso}} (s^{-1})$, still undetermined. The $S_1$ molecule diffuses with the rate $k_{\rm{diff}}  (s^{-1})$, and the molecule $M_2$ occupies the coarsest mesh in the hierarchy consisting of one voxel. This means that $M_2$ does not diffuse, and that $M_1$ always occupies a voxel contained in the voxel of $M_2$.

Now, in each step the probability that the molecules react, $P_{\rm{react}}$, is given by
\begin{equation}
P_{\rm{react}} = \frac{k_r^{\rm{meso}}}{k_r^{\rm{meso}}+k_{\rm{diff}}}.
\end{equation}
Each event has a waiting time, $t_{\rm{event}}$, of
\begin{equation}
t_{\rm{event}} = \frac{1}{k_r^{\rm{meso}}+k_{\rm{diff}}}.
\end{equation}
Thus, the average time $\tau$ until the molecules react is the average number of events until a reaction fires, $(P_{\rm{react}})^{-1}$, times the average time per event, $t_{\rm{event}}$, so we get
\begin{equation}
\tau = (P_{\rm{react}})^{-1}t_{\rm{event}} = \frac{1}{k_r^{\rm{meso}}}.
\end{equation}
We now want to choose $k_r^{\rm{meso}}$ so that we obtain the correct mean binding time $\tau$. However, for a large enough domain, the mean binding time is $(k_{\rm{meso}})^{-1}$, where $k_{\rm{meso}}$ is the reaction rate for the voxel occupied by the molecule $M_2$. Thus the reaction rate should be chosen as if both molecules occupied the coarser mesh, in order to reproduce the correact mean binding time of the molecules.

This argument holds in general. Consider the case where the $M_2$ molecule occupies a level in the hierarchy that is not the coarsest. Assume that the $M_2$ molecule diffuses with diffusion rate $D_A$, and that the $M_1$ molecule is fixed inside a voxel. Now there exists only one voxel that the $M_2$ molecule can occupy, that is also occupied by the $M_1$ molecule. This means that the molecules should react as if both molecules occupy the mesh of the $M_2$ molecule.

\subsubsection{Choosing an initial mesh hierarchy}
\label{choosehier}

It has been shown that mesoscopic simulations of a reversible reaction become more accurate down to some mesh size $h^{\ast} \approx \frac{2}{3}\pi C_3\sigma \approx 3.2\sigma$, where $\sigma$ is the sum of the reaction radii of a reactive pair of molecules \cite{hellanderReactionRatesMesoscopic2015}. For mesh sizes below $h^{\ast}$, simulations actually get less and less accurate, so the optimal mesh size is $h^{\ast}$. It was also shown that for this mesh size, we will reproduce the correct average rebind time.

The relative error of the mean rebind time, $\mesoreact$, for two particles in the RDME using mesh size $h$ is given by \cite{hellanderMesoscopicmicroscopicSpatialStochastic2017}

\begin{align}
\label{eq:W}
W(h) &= \frac{\left| \mesoreact-\microreact \right|}{\microreact} =  \frac{\krmi}{D}G(h,\rrad),
\end{align}
where
\begin{align}
G(h,\sigma) = \begin{cases}
\frac{1}{2\pi}\log\left(\pi^{-\frac{1}{2}}\frac{h}{\sigma}\right)-\frac{1}{4}\left(\frac{3}{2\pi}+C_2\right)\,\,\mathrm{(2D)}\\
\frac{1}{4\pi\sigma}-\frac{C_3}{6h}\,\,\mathrm{(3D)}
\end{cases}
\end{align}
and
\begin{align}
C_d \approx \begin{cases}
0.1951,\,\, d=2\\
1.5164,\,\, d=3
\end{cases}
\end{align}

This error was used to, given a fixed mesh with mesh size $h$ and a chemical reaction system, partition a model into a mesoscopic and microscopic subset in a hybrid method \cite{hellanderMesoscopicmicroscopicSpatialStochastic2017}. Here, we can instead use it to, given a model, compute a largest $h$ for which any given reaction can be handled to satisfy an error $W(h)<\epsilon$, where $\epsilon$ is a user supplied tolerance. This holds when \cite{hellanderMesoscopicmicroscopicSpatialStochastic2017}   

\begin{align}
\label{tolineq}
\krmi (1+\epsilon)^{-1} < \krme h^3
\end{align}
where $\krmi$ is the microscopic reaction rate, and $\krme$ is the mesoscopic reaction rate.

For each species we can compute the coarsest mesh resolution satisfying \eqref{tolineq}. Whenever a molecule is created it will be initialized to this mesh size, and this will be the finest mesh size on which we will ever need to simulate this molecule.

\section{Numerical Experiments}

In this section we show that for some problems we obtain simulation results as accurate as with a microscale simulation, but with a three orders of magnitude speed-up compared to microscale simulations. Specifically we consider problems where molecules can participate in reactions soon after a dissociation, requiring detailed spatial resolution in order to resolve the spatial correlation of the products following such a dissociation.

All timing data has been generated on a Macbook Pro 2017, 3.1 GHz CPU with 8GB of RAM. The problems have been chosen such that \eqref{tolineq} is satisfied to a  sufficiently small $\epsilon$ only for the finest possible mesh size $h^\ast$. All results for the eGFRD algorithm have been generated with a state-of-the-art implementation \cite{PhysRevLett.97.230602, vanzonGreenSfunctionReaction2005, sokoegfrd}, available for download at \url{https://github.com/gfrd/modern_egfrd}.

\subsection{Rebinding dynamics}
\label{sec:rebindingdyn}

In the first example we show that we accurately reproduce the rebinding dynamics of a bimolecular reaction. We thus consider the rebind dynamics of the simple system
\begin{align}
\label{first-system}
\ce{S1 + S2 <-->[k1][k2] S3}.
\end{align} 
Depending on the association rate $k_1$, the products $S_1$ and $S_2$ may have a high probability of rebinding fast. On the microscopic scale, the molecules are placed in contact following a dissociation event, potentially leading to many fast rebind events. On a coarse mesoscopic mesh, the products are placed in the same voxel, but we assume that they are immediately well-mixed inside that voxel. We thus lose spatial information, and will see fewer fast rebinding events.

With the hRDME, we want to accurately reproduce the behavior of a simulation on the finest mesh size $h^*$. This means that the average rebind time will be correct, that the distribution of rebinding times should match the distribution of a mesoscopic simulation on the mesh size $h^*$ (but not necessarily the rebinding-time distribution of a microscopic simulation on length scales smaller than $h^{\ast}$).

\begin{figure}
	\includegraphics[width=0.85\textwidth]{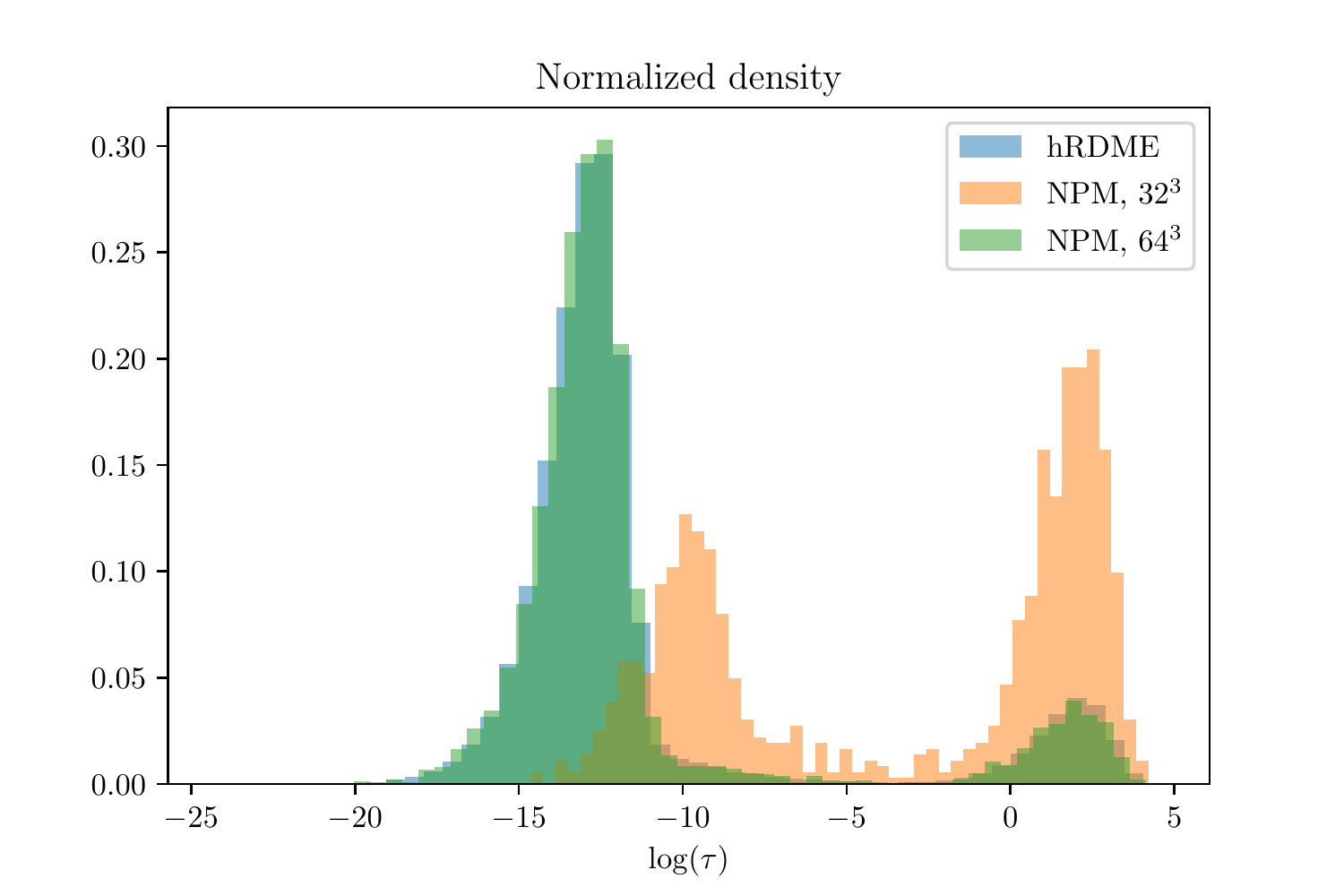}
	\caption{\label{fig-rebind1}The distribution obtained with the hRDME on a sequence of seven meshes overlaps the distribution obtained with the RDME on a mesh of $64^3$ voxels, the finest mesh used for the hRDME simulation.}
\end{figure}

In Fig. \ref{fig-rebind1} we show that for appropriately chosen method parameters, simulations with the hRDME on a sequence of seven meshes, $1$, $2^3$, $4^3$, $8^3$, $16^3$, $32^3$, $64^3$ voxels, is able to reproduce the distribution of rebinding times obtained with a pure RDME simulation on a mesh of $64^3$ voxels. 

The reaction radius of all species is $\sigma = 0.00246$ (so that $h^\ast\approx 3.2\cdot 2\sigma \approx 1/64$), the diffusion constant is $D = 1.0$, and the association rate $k_2 = 1.0$, with a domain volume of 1. Note that the dissociation rate $k_1$ is not important for the rebinding time distribution.

\begin{figure}
	\includegraphics[width=0.85\textwidth]{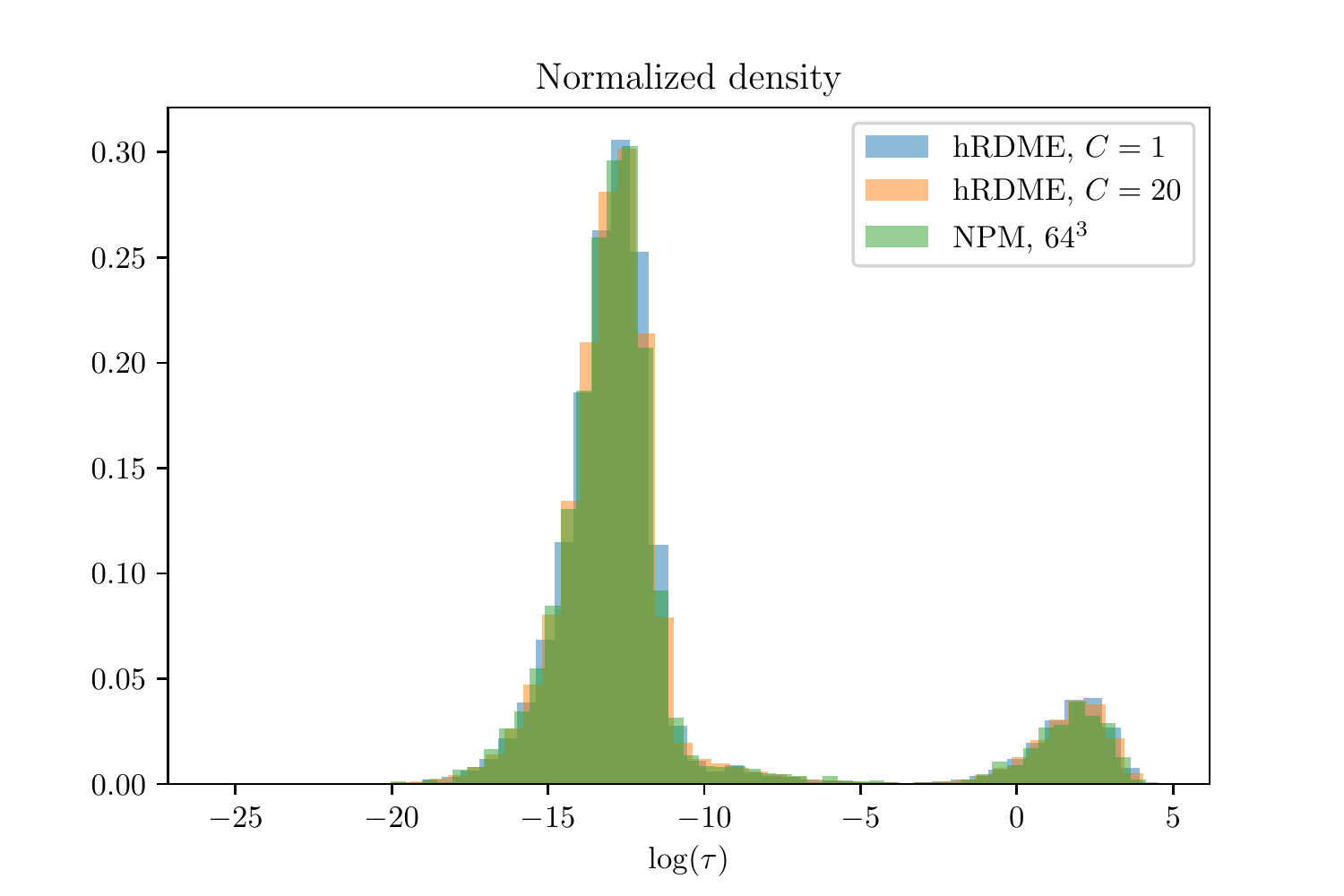}
	\caption{\label{fig-rebind2}Rebinding time distributions for the hRDME with $C=1$ and $C=20$, compared to a pure NPM simulation on a maximally resolved mesh of $64^3$ voxels. The distributions overlap. While $C=1$ could provide sufficient accuracy for many problems, we choose $C=20$ to ensure a large enough $C$ for the vast majority of problems.}
\end{figure}

\paragraph{How to choose the constant C?} In \eqref{Ceq} there is a constant $C$ that controls how much the molecules should diffuse (on average) before they are moved between meshes. To reproduce the rebind distribution we find that $C=1$ seems sufficient. We show this in Fig. \ref{fig-rebind2}. However, as we have no method to determine the optimal value for $C$ for the general case, we choose $C=20$ for the following numerical examples. While we could likely choose a smaller $C$, and thus save even more computational time, we want to choose a $C$ that is likely to work for almost any system.

\subsection{Fast rebinding}

Consider a simple system
\begin{align}
\label{rebind-one}
\ce{S_1 ->[k_1] S_{11} + S_{12} ->[k_2] S_2}.
\end{align}

The same system, and extensions of it, have been studied in detail before \cite{hellanderReactionRatesMesoscopic2015,hellanderMesoscopicmicroscopicSpatialStochastic2017}. Capturing the mean behavior of this system relies heavily on being able to capture the reaction dynamics of $S_{11}$ and $S_{12}$ to sufficient accuracy. To do this, the system has to be simulated on a sufficiently fine mesh. The required resolution can be determined from the criteria given by (\ref{eq:W}) for some sufficiently small $\epsilon$. We showed that $\epsilon < 0.025$ is a reasonable choice for this, and similar, systems \cite{hellanderMesoscopicmicroscopicSpatialStochastic2017}.

If the association reaction is fast, we will need a mesoscopic mesh of maximum resolution. This will make the simulation very expensive, often more expensive than a simulation with the eGFRD algorithm. However, by using the fact that most of the system can be simulated at a coarse level we can speed up the simulation by several orders of magnitude.

For the simple system in \eqref{rebind-one} we note that both the molecules $S_1$ and $S_2$ can be safely simulated on the coarsest scale, so all molecules of these species will be initialized on the coarsest mesh. When an $S_1$ molecule dissociates, the products $S_{11}$ and $S_{12}$ are placed in the same voxel on the finest mesh. We thus resolve the possible rebind events to the highest possible accuracy, and then, if the molecules survive for some time, we can start moving them up in the hierarchy.

All molecules have a reaction radius $\sigma = 0.0025$ and diffuse with diffusion constant $D = 1.0$. We let $k_1 = 1.0$ and $k_2 = 1.0$. The coarsest mesh has 1 voxel, with 7 meshes total in the hierarchy, so that the finest mesh has $(2^6)^3 = 64^3 = 262144$ voxels. The total volume of the domain is $V = 1$, and we simulate the system for a total of $\SI{5}{s}$ and sample the time series in 100 equidistant points between 0 and 5.

We compute the error $E$ as the sum of the difference between the time series, where the eGFRD simulations are considered the correct solution. 

The  hRDME simulation is roughly 2000 times faster than eGFRD, and 500 times faster than a pure mesoscopic simulation of maximum resolution. The accuracy is comparable to the accuracy of the eGFRD simulation. In Table \ref{rebind-table} we present the accuracy of pure NPM simulations on different mesh sizes compared to the hRDME and eGFRD, as well as the wall time per trajectory. The hRDME outperforms both eGFRD and highly resolved mesoscopic simulations, without losing too much accuracy.

A simple extension of the system \eqref{rebind-one} is to add another layer,
\begin{align}
\label{rebind-double}
\ce{$S_1$ ->[$k^1_1$] S_{11} + S_{12} ->[$k^1_2$] S_2}\\
\ce{$S_2$ ->[$k^2_1$] S_{21} + S_{22} ->[$k^2_2$] S_3}
\end{align}

Similarly, products produced by a dissociation is placed in the same voxel on the finest mesh, while $S_1$, $S_2$ and $S_3$ can be safely simulated on the well-mixed scale. This means that we lose spatial information in between reactions, but this will not negatively affect the accuracy as long as we accurately capture fast rebinding events following dissociations. All parameters are as above, with all reaction rates equal to 1.0.

\begin{figure}
	\includegraphics[width=0.49\textwidth]{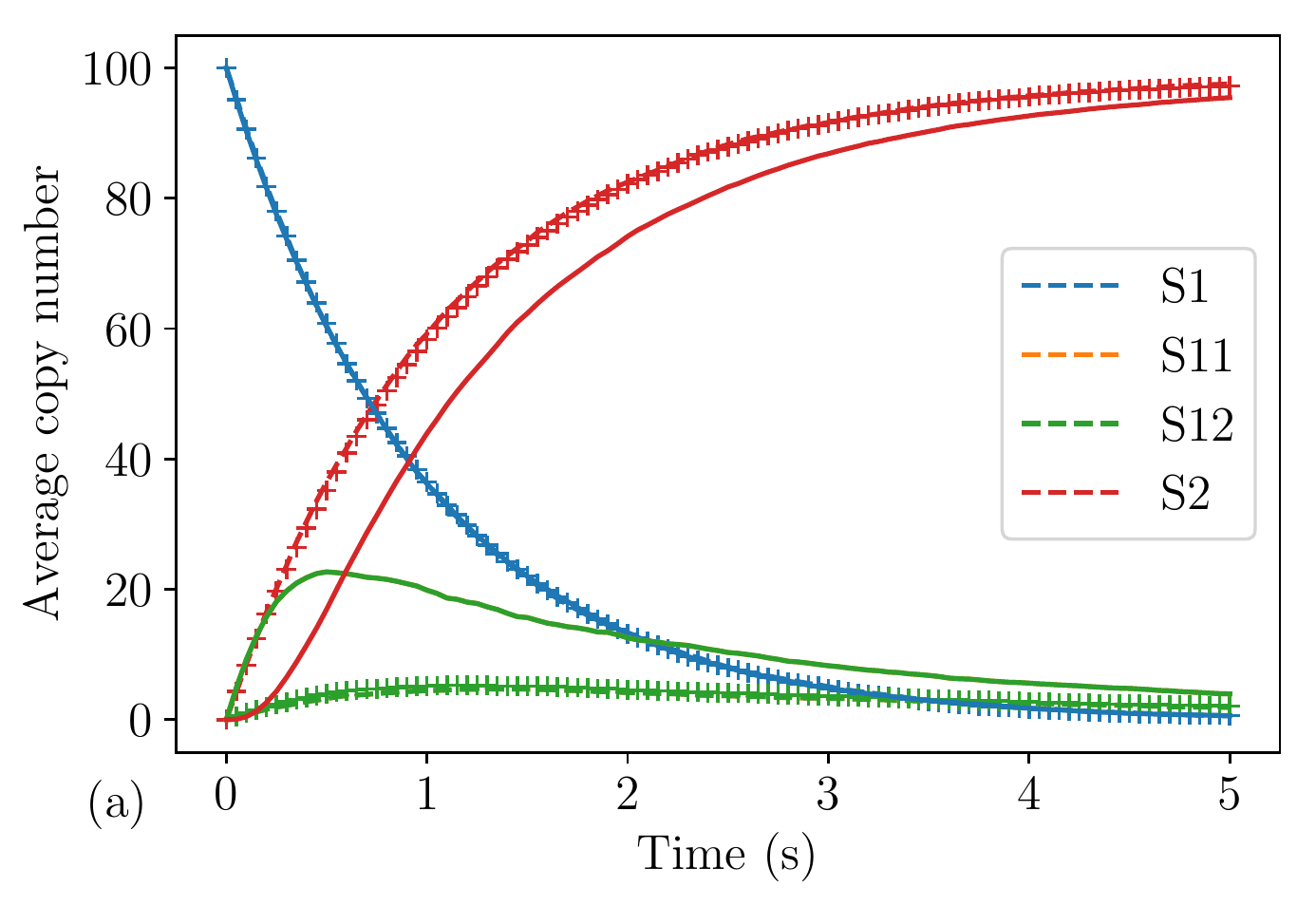}
	\includegraphics[width=0.49\textwidth]{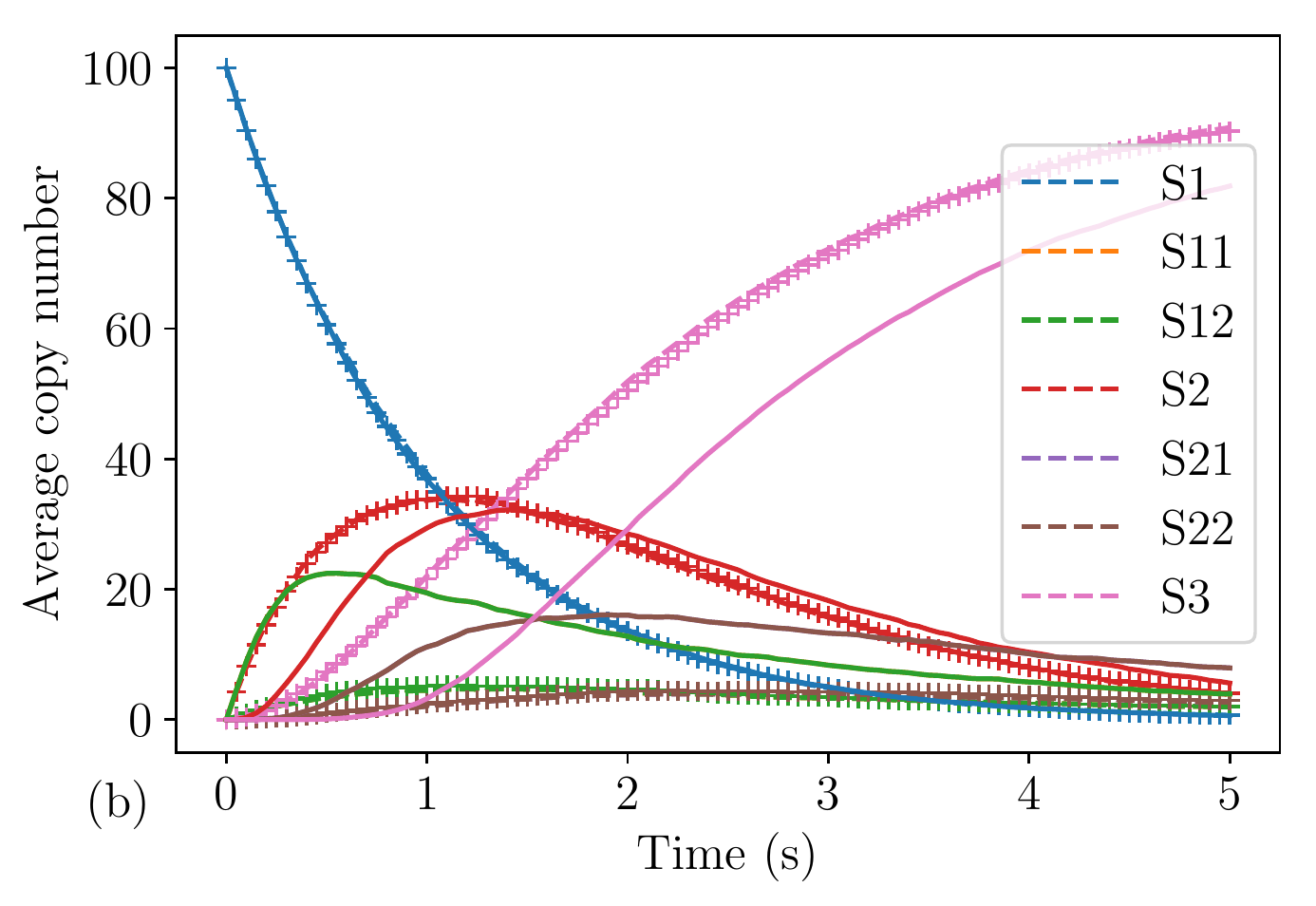}
	\caption{\label{fig-rebind-double}{In (a) we plot the time series of the system \eqref{rebind-one}, and in (b) the time series of the system \eqref{rebind-double}. The eGFRD results are plotted with dashed lines, and hRDME with $(+)$. For reference we have plotted the well-mixed results in solid lines. As we can see, eGFRD results and hRDME results agree well, while there is a significant error in the well-mixed simulation results. The average is based on 200 trajectories, sampled at 101 points from 0 to 5.}}
\end{figure}

The speed-up compared to eGFRD is in this case approximately a factor of 1800, and 400 times faster than an NPM simulation on a maximally resolved mesh. In Fig. \ref{fig-rebind-double} we plot the time series of both of the systems above, simulated with the hRDME and eGFRD. There is no visible difference between the results. In Table \ref{rebind-table} we present the error (computed as above) and wall time per trajectory for simulations with the NPM on different mesh sizes, compared to the hRDME and eGFRD.

\begin{table}
	\setlength{\tabcolsep}{10pt}
	\begin{tabular}{l|lcccccc}
		\hline
		{} System &  &  $8^3$ &  $16^3$ &  $32^3$ &  $64^3$ &  hRDME &   eGFRD \\
		\hline
		\multirow{2}{*}{Single  \eqref{rebind-one}} & Speedup  & 8.45 & 30.60 & 134.35 & 528.33 & 1.00 & 2006 \\
		& Error     & 0.1529 &  0.1292 &  0.0936 &  0.0031 & 0.0062 &   \\
		\hline
		\multirow{2}{*}{Double  \eqref{rebind-double}} & Speedup &         7.46 &        26.39 &       110.93 &       404.95 &   1.00 & 1800 \\
		& Error & 0.2748 &  0.2419 &  0.1684 &  0.0079 & 0.0072 &   \\\hline
	\end{tabular}
	\caption{Speedup, as a multiple of the hRDME (so that a large number means a slower simulation, and a small number is faster), and relative error. We have tabulated the results of simulations with the NPM on a single mesh with varying mesh resolution ($n^3$ corresponds to a simulation with the NPM on a Cartesian mesh consisting of $n^3$ voxels), results of the hRDME and finally results from simulations with the eGFRD algorithm. The estimate of the error is based on 200 trajectories.}
	\label{rebind-table}
\end{table}

\subsection{MAPK}

Takahashi et al. have shown that a MAPK system, for some parameter values, exhibits a fine-grained dynamics that cannot be accurately resolved with a well-mixed model \cite{TaTNWo10}. Hellander et al. were able to reproduce the behavior of the system with a highly resolved RDME simulation \cite{hellanderReactionRatesMesoscopic2015}. However, this required a maximally resolved mesh, in that case $64^3$ voxels, making the simulation very slow and completely dominated by diffusion events. At this resolution, the RDME was slower than microscale simulations.

The MAPK model is given by
\begin{align}
\label{MAPK-system}
\ce{KK + K <-->[k_1][k_2] KK&-K ->[k_3] KK^{\ast} + K_p}\\
\ce{KK + K_p <-->[k_4][k_5] KK&-K_p ->[k_6] KK^{\ast} + K_{pp}}\\
\ce{P + K_{pp} <-->[k_1][k_2] P&-K_{pp} ->[k_3] P^{\ast} + K_{p}}\\
\ce{P + K_{p} <-->[k_4][k_5] P&-K_p ->[k_6] P^{\ast} + K}\\
\ce{KK^{\ast} &->[k_7] KK}\\
\ce{P^{\ast} &->[k_7] P}.
\end{align}

The volume of the domain is $V = 1.0$, all species diffuse with diffusion constant $D=1.0$ and the reaction radius of all species is $\sigma = 0.0024599$. The reaction rates are
\begin{align}
k_1 &= 0.0448346\\
k_2 &= 1.35\\
k_3 &= 1.5\\
k_4 &= 0.0929902\\
k_5 &= 1.73\\
k_6 &= 15.0\\
k_7 &= 693147.18.
\end{align}

Here we show that it is possible to simulate this system with high accuracy, but with a great speed-up compared to GFRD simulations. For simplicity, we have made no assumptions about which species could be simulated on a coarse-grained level. It is possible that the simulations could be optimized even more by initializing some molecules on a mesh that is coarser than the finest mesh.

\begin{figure}
	\includegraphics[width=0.85\textwidth]{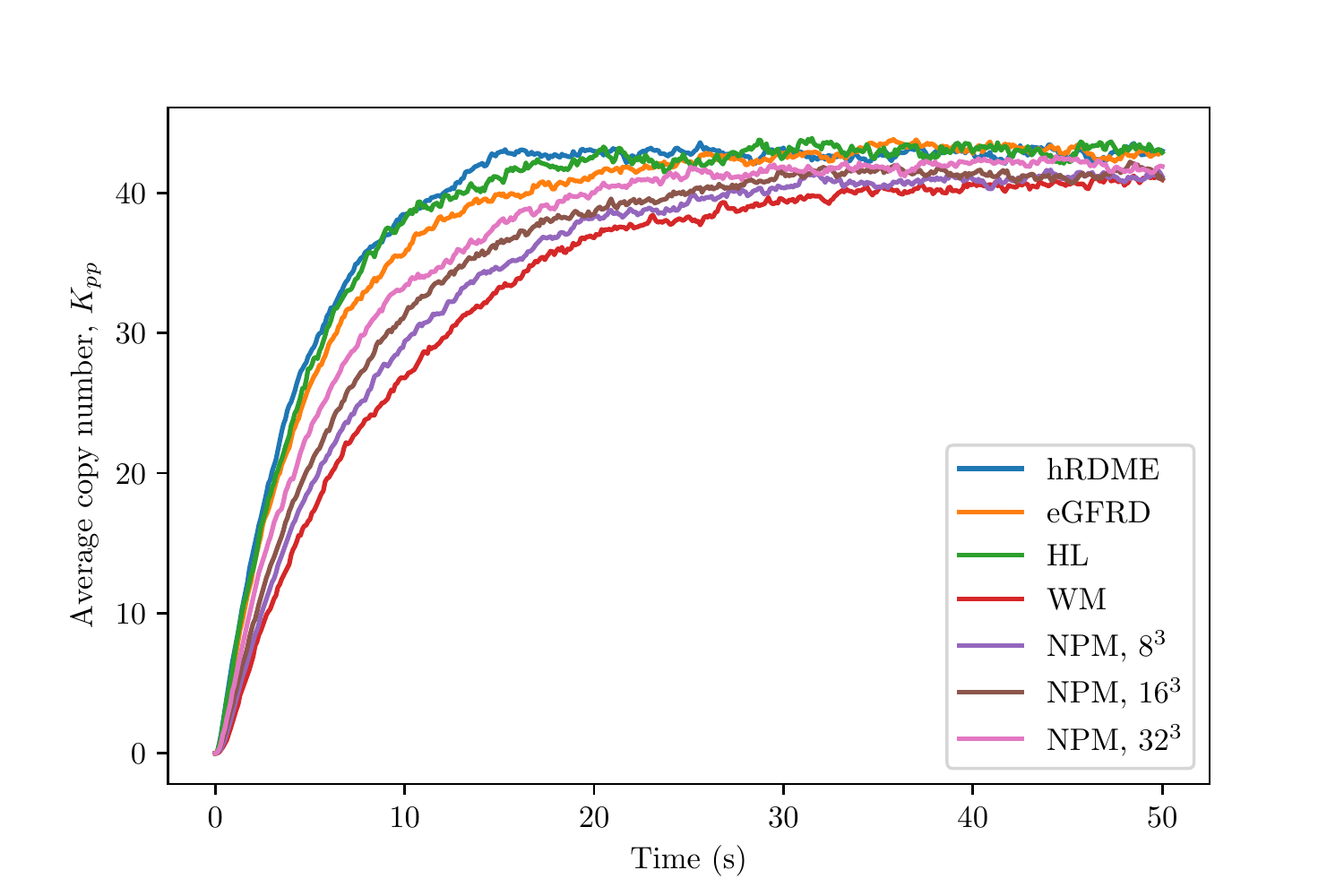}
	\caption{\label{fig-mapk-ts}{The system was simulated for \SI{50}{s}. The copy numbers are the average of 200 trajectories. We simulated the system with the hRDME algorithm, an implementation of the eGFRD algorithm, and with a microscale solver implemented by the authors \cite{hellanderSingleMoleculeSimulations2013} (based on the same modeling framework as eGFRD, and denoted by HL in the plot above). As we can see, all simulations match reasonably well, but there is a small difference compared to the eGFRD results, while the hRDME  matches the other microscale implementation (HL) very well. We have also plotted results of pure mesoscopic simulations on different mesh sizes for reference.}}
\end{figure}

In Fig. \ref{fig-mapk-ts} we plot the average time series of the species $K_{pp}$. The hRDME simulations matches the microscale simulations well. We have also simulated the system with the NPM on different mesh resolutions, for reference. Timing results are presented in Table \ref{tab:mapk}.

\begin{table}
	\setlength{\tabcolsep}{10pt}
	\begin{tabular}{ l c c c c c c c }
		\hline
		& WM &  $8^3$ &  $16^3$ &  $32^3$ & $64^3$ &  hRDME &    eGFRD \\
		\hline
		Speedup & 0.05 &        0.35 &         1.36 &         5.36 &        94.15 &   1.00 &  59.03 \\
		\hline
	\end{tabular}
	\caption{Speedup (defined as in Table \ref{rebind-table}) for different mesh resolutions for the MAPK system, compared to eGFRD simulations. The speed-up is roughly a factor of 60 compared to eGFRD, with a wall time of \SI{7.51}{s} per trajectory for hRDME. We have also simulated the system with $C=1$ in the hRDME algorithm, with a total speed-up of 170 and no noticeable difference in accuracy compared to the simulations with $C=20$.}
	\label{tab:mapk}
\end{table}

\section{Discussion}

We have shown that coupling mesoscopic simulations on different mesh sizes can save orders of magnitudes of computational time, while being as accurate as microscale or highly resolved mesoscopic simulations. This methodology is also faster, and much simpler to implement, than a mesoscopic-microscopic hybrid scheme. However, there are still cases where a hierarchical mesoscopic simulation will not be sufficiently accurate. In particular, if molecular crowding effects are important, they are not captuted by the mesoscopic model, while they are captured in the microscopic hard-sphere model. The accuracy is also still limited by $h^*$ \cite{hellanderReactiondiffusionMasterEquation2012,hellanderReactionRatesMesoscopic2015}; this lower bound on the mesh size is inherent to the mesoscopic model and not due to the methodology presented herein.

Here we have considered structured Cartesian meshes only. While the methodology could in principle be extended to unstructured meshes, this is technically more difficult, and the complex shape of the voxels would incur a larger overhead. This in turn means that it is more difficult to handle complex geometries than it is when we have pure microscopic or mesoscopic simulations, or hybrid methods.

\begin{acknowledgments}
The work has been funded by the Swedish Research Council, the eSSENCE strategic collaboration on eScience, and the NIBIB of the NIH under grant no. NIH/2R01EB014877-04A1. 
\end{acknowledgments}

\bibliography{hrdme}

\end{document}